\documentclass[12pt]{amsart}
\usepackage{amsthm,amsmath,amssymb,amscd,graphics,enumerate, ifthen}

\def\version{SHORT}

\newcommand{\marg}[1]{}
\newcommand{\thisdate}{\today}
{
   \newtheorem{theorem}[subsubsection]{Theorem}
   \newtheorem{proposition}[subsubsection]{Proposition}     
   \newtheorem{lemma}[subsubsection]{Lemma}
   \newtheorem*{claim}{Claim}

}
{\theoremstyle{definition}

   \newtheorem{definition}[subsubsection]{Definition}
   
}
\newcommand{\RR}{{\mathbb{R}}}

\newcommand{\CC}{{\mathbb{C}}}
\newcommand{\QQ}{{\mathbb{Q}}}

\newcommand{\PP}{{\mathbb{P}}}
\newcommand{\ZZ}{{\mathbb{Z}}}

\newcommand{\GG}{{\mathbb{G}}}
\newcommand{\LL}{{\mathbb{L}}}
\newcommand{\bbA}{{\mathbb{A}}}

\newcommand{\bK}{{\mathbf{K}}}
\newcommand{\bM}{{\mathbf{M}}}

\newcommand{\bX}{{\mathbf{X}}}

\newcommand{\bmu}{{\boldsymbol{\mu}}}

\newcommand{\cB}{{\mathcal B}}
\newcommand{\cC}{{\mathcal C}}

\newcommand{\cG}{{\mathcal G}}

\newcommand{\cK}{{\mathcal K}}

\newcommand{\cM}{{\mathcal M}}

\newcommand{\cO}{{\mathcal O}}

\newcommand{\cS}{{\mathcal S}}

\newcommand{\cX}{{\mathcal X}}
\newcommand{\cY}{{\mathcal Y}}

\newcommand{\Spec}{\operatorname{Spec}}

\newcommand{\Pic}{{\operatorname{\mathbf{Pic}}}}
\newcommand{\Hom}{{\operatorname{Hom}}}

\newcommand{\Ext}{{\operatorname{Ext}}}

\newcommand{\Aut}{{\operatorname{Aut}}}

\newcommand{\lrar}{\longrightarrow}
\newcommand{\mto}{\mathop\to\limits}
\newcommand{\dar}{\downarrow}

\newcommand{\ocM}{\overline{{\mathcal M}}}
\newcommand{\obM}{\overline{{\mathbf M}}}

\newcommand{\ocX}{\overline{{\mathcal X}}}

\newcommand{\e}{{\mathbf{e}}}

\newcommand{\double}{\genfrac..{0pt}1
{\raise -1pt\hbox{$\scriptstyle\longrightarrow$}}{\raise 3pt\hbox
{$\scriptstyle\longrightarrow$}}} 
\newcommand{\setmin}{\,\protect%
\begin{picture}(8,3.5)\qbezier(1,3.5)(4,2.)(7,.5)\end{picture}\,}

\def\gen{_{\rm gen}}

\def\tototi{\mathbin{\mathop{\otimes}\limits^{\raise-1pt\hbox
{$\scriptscriptstyle {\rm L}$}}}}

\def\indlim{\mathop{\vrule width0pt height7pt depth
4pt\smash{\lim\limits_{\raise 1pt\hbox to 14.5pt
{\rightarrowfill}}}}}
\def\projlim{\mathop{\vrule width0pt height7pt depth
4pt\smash{\lim\limits_{\raise 1pt\hbox to 14.5pt
{\leftarrowfill}}}}}

\begin{document}
\title[Quantum products]{Algebraic orbifold quantum products}
\author[D. Abramovich]{Dan Abramovich}
\thanks{Research of D.A. partially supported by NSF grant DMS-0070970}  
\address{Department of Mathematics\\ Boston University\\ 111 Cummington
         Street\\ Boston, MA 02215\\ U.S.A.} 
\email{abrmovic@math.bu.edu}
\author[T. Graber]{Tom Graber}
\thanks{Research of T.G. partially supported by an NSF post-doctoral
research fellowship. }
\address{Department of Mathematics\\ 
Harvard University\\ 
1 Oxford Street\\
Cambridge MA 02138\\ 
U.S.A.} 
\email{graber@math.harvard.edu}
\author[A. Vistoli]{Angelo Vistoli}
\thanks{Research of A.V.  partially supported by the University of
        Bologna, funds for selected research topics.}
\address{Dipartimento di Matematica\\ Universit\`a di Bologna\\Piazza di Porta
         San Donato 5\\ 40127 Bologna\\ Italy}
\email{vistoli@dm.unibo.it}
\date{\thisdate}

\maketitle

\section{Introduction}

The purpose of this note is to give an overview of our work on defining algebraic
counterparts for W.\ Chen and Y.\ Ruan's Gromov-Witten Theory of
orbifolds. This work  will be described in detail in a subsequent
paper. 

The presentation here is generally based on lectures given by two of
us at the Orbifold Workshop in Madison, Wisconsin. Following the
spirit of the workshop, our presentation is intended to be
understandable not only to Algebraic geometers, but also practitioners
of the other disciplines represented there, including differential
geometers and mathematical physicists.

We make a special effort to make our constructions as canonical as
we can, systematically using the language of algebraic stacks. Our
constructions are based on the theory of twisted stable maps developed
in \cite{A-V:stable-maps}, but requires making explicit some details which
were not studied in the paper \cite{A-V:stable-maps}. Apart from the
pleasure we take in understanding these details, our efforts bear some
concrete fruits; in particular, we are able to define the Chen-Ruan
stringy product in degree 0 (the so called stringy cohomology) with
integer coefficients.


 We work over the field $\CC$ of complex numbers (although our
 discussion works just as well over any field of characteristic 0).

\section{Stacks and their moduli spaces}
 
There are two, quite different, ways in which ``orbifolds'' or
``stacks'' arise.

\subsection{Groupoids in schemes} This is the way most differential
geometers, as well as  many algebraic geometers, are introduced to the
subject, since it is in some sense concrete and  geometric: one thinks
about an object which is locally modeled on a quotient of a variety by
the action of an algebraic group. Then one needs to define a good
notion of maps between such objects - this is the difficult part of
the picture.

Concretely, one is given a ``relation'', namely a morphism $R \to U\times U$, 
where each of the projections $R \to U$ is \'etale (or, more generally, smooth), 
with the extra data of an ``inverse'' $R \to R$ and ``composition'' $R
\times_U R \to R$ satisfying natural assumptions which the reader may
guess. This is what one called a groupoid. 

One then defines morphisms of groupoids, Morita equivalences,
and finally one defines the category of orbifolds using groupoids up
to Morita equivalence.  

In the theory of stacks, there is a natural algebraic stack
$[U/R]$  associated to a groupoid. Two groupoids are Morita
equivalent if and only if the associated stacks are isomorphic; therefore
the groupoid appears as a \emph{presentation} of the stack. The stack is a
more intrinsic object.

For instance, the quotient stack $[U / G]$ of a scheme $U$ by the
action of a finite group $G$ is the stack associated to the groupoid
given by  $R = G \times U$.

A particularly simple and useful case is the {\em classifying stack}
$\cB G = [pt / G]$ of the group $G$ (where $pt$ stands for the point
$\Spec \CC$).

\subsection{Moduli problems} This is the way stacks arise in a more
abstract, categorical incarnation, but still extremely natural.

One can encode the data of a ``moduli problem'' in the {\em category}
 of families of objects which we want to parametrize. For instance,
 the data of ``moduli of curves of genus $g$''  is encoded in the
 category $\cM_g$ whose objects are 
  $$\cO b (\cM_g) = \left\{\ \ \begin{array}{c} C\\ \dar \\
S\end{array} \ \ \left\bracevert \hbox{ 
\parbox{4in}{\begin{center}$S$ is a 
scheme, \\ and \\ $C \to S$ is a smooth family of curves of genus
$g$\end{center}}}\right\}\right.,$$
and the arrows from one object $C \to S$  to another $C' \to S'$ are
the {\em fiber diagrams} 
$$ \begin{array}{ccc} C & \lrar & C' \\ \dar & & \dar \\ S & \lrar& S'.\end{array}
$$
There is also a {\em structural functor} $\cM_g \to Sch$ which sends a
family of curves $C \to S$ to the underlying scheme $S$.

 Formally, an
 {\em algebraic stack} 
 is a category $\cM$ with a functor  $ \cM \to Sch$, satisfying
 certain assumptions which guarantee that it is ``somewhat close to being
 representable by a scheme''. These assumptions imply, in particular,
 that an algebraic stack is always the stack associated to a groupoid
 in schemes $[U/R]$. If the automorphism group of every object in the
 category is {\em finite}, then in fact  $\cM$  is {\em locally} given as the
 quotient $[U/G]$ of a scheme by a {\em finite group}. Such a stack is
 called a {\em Deligne--Mumford stack}.
 
 If $\cM$ is a stack of a moduli problem, then exactly as in the example, the
 objects are families $C \to  S$, and  the arrows are fiber diagrams
 $$\begin{array}{ccc} C & \lrar & C'\\ 
\dar &&\dar \\   B & \lrar & S'. \end{array}$$ The functor $\cM \to
 Sch$ simply sends a family $C \to S$ to its base scheme $S$.

 Simple examples include: $\cM_g,$ the moduli stack of smooth curves of
 genus $g$, and  $\cB G$, the moduli stack of principal homogeneous
 $G$-spaces, say for a finite group $G$.

In introductory texts about moduli theory, one learns about the
{\em moduli functor} of a moduli problem  - in the example on $\cM_g$ above, the
functor sends a scheme $S$ to the {\em set of isomorphism classes} of
all smooth families $C \to S$ of curves of genus $g$. It may be the
result of an 
unfortunate historical tradition, that often the discussion of the {\em
category}, which faithfully encodes the moduli problem in question, is
delayed till after the shortcomings of the moduli functor are discovered.

\subsection{Coarse moduli spaces}

Every Deligne--Mumford stack $\cM$ has an associated algebraic space
called the {\em coarse moduli space
$\bM$} (for a general proof see \cite{Keel-Mori}):
 
\begin{itemize} 
\item In case $\cM$ is given via a groupoid as $[U/R]$, the coarse moduli space
is the orbit space $U/R$.

\item
In case $\cM$ is the moduli stack of a moduli problem, then 
$\bM$ is the usual coarse moduli space parametrizing  objects modulo
isomorphisms. 
\end{itemize}

\section{Twisted stable maps}
\subsection{Stable maps into a projective variety}
Recall that, when $X$ is a projective variety, one has a {\em
projective} moduli space
$\obM_{g,n}(X,\beta)$ 
of {\em $n$-pointed stable maps of genus $g$ and image class
$\beta$.} This is always  the coarse moduli space of a
Deligne--Mumford stack $\ocM_{g,n}(X,\beta)$,
whose objects over a scheme $T$ are triples:
$$ \left(C \to T, \ \ \{s_1, \ldots s_n: T \to C\}, \ \ f: C \to
X\right). $$
Here \begin{enumerate}
\item $C \to T$ is a family of prestable curves of arithmetic genus $g$, 
\item $s_i: T \to C$ are sections of $C \to T$ with disjoint images
lying in the smooth locus, and
\item $f: C \to X$ is a morphism, such that the group of automorphism
of fibers of $C\to T$ commuting with $f$ and fixing all the sections
$s_i$ is finite.
\end{enumerate}

It is no secret in this workshop that the spaces (and stacks) of
stable maps have been useful in symplectic geometry,  enumerative
geometry and mathematical 
physics through Gromov--Witten  invariants. They have also served as a
stepping stone for constructing objects of interest in algebraic
geometry (such as other moduli spaces). For a recent application in
higher dimensional geometry see \cite{Graber-Harris-Starr}.   

\subsection{Stable maps into a stack}
Since here we are in the business of enlarging our vocabulary of geometric
objects from ``varieties'' to ``orbifolds'' (or ``stacks''), one may
wonder if it is worthwhile to replace the projective variety $X$
in the definition of a stable map by a stack $\cX$. 

Let us consider two examples of stacks which were mentioned above: 
\begin{enumerate} \item Fix an integer $h>1$, and say $\cM =
\ocM_h$. Then a map of a curve $C$ to $\cM$ is the same as an object
of $\cM$ over $C$, in other words a {\em   family $X \to C$ of stable curves
of genus $h$ over $C$}. Thus the moduli space of all such things is in
particular a moduli space of a certain type of surfaces with extra
structure (given by a map to a curve). For an algebraic geometer,
it is evident that moduli of surfaces are of interest. So having a
space of stable maps into $\ocM_g$ is useful at least as a
construction tool. See \cite{A-V:fibered-surfaces}.
\item Fix a finite group $G$, and say $\cM = \cB G$. Then a map of a
curve $C$ to $\cM$ is a family $P\to C$ of principal homogeneous
$G$-spaces, in other words, a principal $G$-bundle. The moduli of
these objects are closely related to  Hurwitz
schemes, level structures  and admissible covers, which have been a
subject of interest in algebraic geometry for years. See \cite{A-C-V}.
\end{enumerate}

It is from this point of view, of using stable maps as a basic
construction tool, that two of us approached the problem of
constructing a suitable space of stable maps into a stack, see
\cite{A-V:stable-maps}. Admittedly, Gromov--Witten theory for stacks seemed to us a
distant possibility until the work of W.\ Chen and Y.\ Ruan on the
subject was made
public.

\subsection{The non-compactness problem} One still needs to give a
definition of what one means by a ``stable map into a stack''. As a
first attempt, one can 
define stable maps into a stack $\cX$ just as stable maps into a
projective variety were defined. The only problem with this is, that
the stack of usual stable maps into $\cX$ is not necessarily compact.

This phenomenon is already manifest with the simplest series of
examples where $\cX = \cB G$, so that maps to $\cX$ correspond to
principal $G$-bundles, see \cite{A-V:stable-maps}. Say $G =
(\ZZ/m\ZZ)^4, n\geq 2,$ and  
$g=2$. A smooth curve $C$ of genus $2$ carries a {\em connected} principal
$G$-bundle $P \to C$, since its first homology group is $\ZZ^4$. But if one
takes a one-parameter family $P_t \to C_t, t\neq 0$ of such bundles
where the underlying curve $C_t$
degenerates, as $t$ approaches $0$, to a nodal curve $C_0$ of geometric genus
1, then there is nothing 
to which the one-parameter family of bundles can degenerate. This is so simply
because the first homology of the limiting curve $C_0$ is $\ZZ^3$, and
thus $C_0$
carries no connected principal $G$-bundle.

An example of a similar phenomenon with $\cX = \ocM_{1,1}$ is described in
\cite{A-V:fibered-surfaces}.

\subsection{Adding orbispace structure on curves} A simple analysis, detailed in
\cite{A-V:families}, \cite{A-V:stable-maps}, of what goes wrong in
the example above, reveals that the problem is centered above the
node of the degenerate curve $C_0$, and, moreover, at the node there
is a natural orbispace structure $\cC_0$ on $C_0$ which does afford a
limiting connected principal bundle $P_0 \to \cC_0$.

One then realizes that, in order to have a {\em compact} moduli space
(or stack) of maps into a stack $\cX$, {\em one needs to allow the
source curves to acquire an orbispace structure as well.} This may seem as a
surprise at first, but we have come to believe that this is completely
natural and beautiful. In a way, nature imposes this solution upon us,
so why fight the elements?

\subsection{Twisted curves}
Let us describe in some detail the orbispace curves - which we call
{\em twisted curves}, underlying these new maps. We restrict ourselves
here to {\em balanced} twisted curves, the only type that is relevant
to our discussion in this paper.

\subsubsection{Nodes} First, consider the structure at a node. The
orbispace curve here is 
locally the quotient $[U/\bmu_r]$ of a nodal curve $U = \{xy=f(t)\}$
by the following action of the
cyclic group $\bmu_r$ of $r$-th roots of 1:
$$(x,y) \mapsto (\zeta x,\zeta^{-1}y).$$

The {\em coarse curve} $C$ underlying the orbispace curve $\cC$ is, locally, the
schematic quotient $U/\bmu_r$ defined by $uv = f(t)^r$, where $u =
x^r$ and $v = y^r$.

This kind of orbispace structure by itself would suffice for
solving the compactification problem we discussed above. But we also
know that it is useful to be able to describe the ``boundary'' of
moduli by gluing {\em marked} curves into nodal ones. Reversing this
line of thought, it is useful to understand what kind of structure one
obtains when {\em separating} a node into two marked points on a maked curve.

\subsubsection{Markings} We are led to consider an orbispace structure along a
marking. Here the orbispace curve $C \to T$ is locally the
quotient $[U/\bmu_r]$ of a smooth  curve $U$, with local
coordinate $z$ defining the marking, by the
following action of $\bmu_r$:
$$z \mapsto \zeta z.$$

 The coarse curve $C$ is locally a smooth curve with coordinate $w = z^r$.

The integer $r$ in these local descriptions is called the {\em index} of
the orbifold curve at the point in question.

\subsubsection{Global structure of a marking} Before imposing on the
reader the formal definition of a twisted 
curve, we need to put a word of caution. Having given the {\em local} description
of a family of twisted curves $\cC \to T$ at a marking, one might be
misled to believe that there is automatically {\em globally} a section
$T \to \cC$. {\em This is not the case in general.} This fact is going
to be important later when we consider evaluation maps.  

Take the coarse family of curves $C \to T$, with its section $T
\to C$, having image $\Sigma^C\subset C$. The reduction
$\Sigma^\cC\subset \cC$ of the inverse image of $\Sigma^C$ in the
twisted curve $\cC$ is canonically a {\em gerbe} over $T$ banded by
the group $\bmu_r$.  This means, in particular, that {\em locally}
$\Sigma^\cC$ is isomorphic to the stack-theoretic quotient of
$\Sigma^C$ by the trivial action of $\bmu_r$, but globally there is an
obstruction for such 
a quotient description. The obstruction, which naturally lies in $H^2(T,
\bmu_r)$, comes in our case from $H^1(T, \GG_m) = \Pic(T)$ in  the
following manner: 

Let $L=s^*N$ be the pullback to $T$ of the normal bundle of $\Sigma^C
\subset C$. Then $\Sigma^\cC$ is {\em the moduli stack of $r$-th roots of
$L$,} whose objects over a scheme $S$ are triples
 $$(g:S \to T,\ \ \ M,\ \ \  \sigma: M^{\otimes r} \tilde{\lrar} g^* L),$$
where 
\begin{itemize} 
\item $g:S \to T$ is a morphism, 
\item $M$ is a line bundle on $S$, and
\item $\sigma$ is an isomorphism of line bundles.
\end{itemize}
The group $\bmu_r$ acts on every object here by fiberwise
multiplication on $M$. All automorphisms arise this way. 

It is
clear from this description that $\Sigma^\cC\to T$ has a section if
and only if the line bundle  $L$ has an $r$-th root over $T$.

Incidentally, there is also a simple prescription for obtaining $\cC$
from the coarse curve $C$, at least away from the nodes. We will not
need that description in the present exposition.


\subsubsection{} Finally, here is a formal definition of twisted
curves, as in \cite{A-V:stable-maps}: 

\begin{definition}
A {\em twisted nodal $n$-pointed curve over a scheme $T$} is a  diagram
$$\begin{array}{ccc} \Sigma^{\cC} & \subset & \cC  \\  &\searrow & \dar  \\ &
 & C   
  \\ &  & \dar  \\ &&T  
\end{array}$$
where
\begin{enumerate}
\item $\cC$ is a tame Deligne-Mumford stack, proper over $T$, which
\'etale locally is a nodal curve 
over $T$;
\item $\Sigma^\cC = \cup_{i=1}^n \Sigma_i^\cC$, where $\Sigma_i^{\cC}
\subset  \cC$ are disjoint closed substacks in the 
smooth locus 
of $\cC \to S$;  
\item $\Sigma_i^{\cC} \to T$ are \'etale gerbes; 
\item the morphism $\cC \to C$ exhibit $C$ as the coarse moduli scheme of 
$\cC$; and
\item $\cC \to C$ is an isomorphism over $C\gen$.
\end{enumerate}
\end{definition}

The notation $C\gen$ in the definition above stands for the {\em
generic locus}, namely the complement of the nodes and markings on
$C$. 

One categorical issue we have to mention is the following: the
collection of all stacks naturally forms a 2-category, and therefore
families of twisted curves also form a 2-category. However, the fact
that 
each twisted curve has a dense open subset which is a scheme, can be
shown to imply that the 2-category of twisted curves is equivalent to
a category (so it has a chance of being a stack).

\subsection{Twisted stable maps}

We can now define a {\em twisted stable map} as follows (see \cite{A-V:stable-maps}):

\begin{definition}
Let $\cX$ be a Deligne--Mumford stack, $\cX \to X \subset \PP^N$
   a  projectively embedded coarse moduli scheme.

An $n$-pointed {\bf twisted stable map} $f:\cC \to \cX$ to $\cX$ of genus $g$,
degree $d$, is a  
diagram 
$$\displaystyle \begin{array}{ccc} \cC & \stackrel{f}{\to} & \cX \\
                \dar & & \dar \\
                C &  \stackrel{\overline{f}}{\to} & X
\end{array}$$ where
\begin{enumerate}
\item $\cC$ is a twisted marked curve, with coarse moduli space  $\cC \to C$,
\item  $C\to X$ a stable map  of
geus $g$ and degree $d$, and 
\item $\cC \to \cX$ representable 
\end{enumerate}
\end{definition}

With this definition, we have the following theorem:

\begin{theorem}\label{Th:stablemaps} The category of $n$-pointed
twisted stable maps of genus $g$
and 
degree $d$ to $\cX$ is a proper Deligne--Mumford stack
$\cK_{g,n}(\cX,d)$ admitting a  
projective coarse moduli space $\bK_{g,n}(\cX,d)$.

\end{theorem}

Of course, the stacks $\cK_{g,n}(\cX,d)$ can be further decomposed
using a numerical equivalence class, or a homology class $\beta\in
H_1(X, \ZZ)$ instead of the degree. 

The notation using the letter $\cK$ rather than $\ocM$ was chosen
originally to avoid confusion when inserting something like $\ocM_\gamma$
as an argument in  place  of $\cX$. In this note we will introduce the notation
$\ocM_{g,n}(\cX,d)$ for a  certain variant of $\cK_{g,n}(\cX,d)$
(namely the
stack of twisted stable maps with trivialized markings).

The stacks $\cK_{g,n}(\cX,d)$ can also be decomposed in terms of the
type of the markings. To 
avoid cumbersome notation, we
will not do this explicitly.  Note, however, that the index
$r_i$ of the  $i$-th marking is a locally constant function
$\cK_{g,n}(\cX,d) \to \ZZ$. As usual with locally constant functions,
$r_i$ induces a linear operator on homology, cohomology, and Chow groups. 

Also, as in the case of usual stable maps, there is a stabilization
morphism $\cK_{g,n}(\cX,*)\to \cK_{g,n}(\cY,*)$ associated to a
morphism of stacks $\cX \to \cY$. This is going to be useful in this
paper in the special case where $\cY$ is a point.

Theorem \ref{Th:stablemaps} is closely related to Proposition 2.3.8 in Chen
and Ruan's paper 
\cite{Chen-Ruan}. There are several differences: Chen and Ruan work
with pseudo-holomorphic maps into smooth symplectic orbifolds in
general. In \cite{A-V:stable-maps} we
work only with algebraic Deligne--Mumford stacks. On the other hand we
allow the target stack $\cX$  to 
be singular (and allow tame stacks in positive characteristics). Chen
and Ruan's moduli space is a differential orbifold, whereas we work
systematically with stacks.

Some words on the proof of the theorem in \cite{A-V:stable-maps}:
often, the construction of a 
moduli stack (or space) is a standard process based on the existence
of Hilbert and Quot schemes. Unfortunately, there is no well-developed
machinery of universal constructions analogous to Hilbert and Quot
schemes for stacks.  The proof in  \cite{A-V:stable-maps} constructs
the stack $\cK_{g,n}(\cX,d)$ ``with bare hands'' using deformation
theory and algebraization, following the list of axioms in M.\ Artin's
paper \cite{Artin}.
 Some of the methods developed in 
\cite{A-V:stable-maps} might serve as building block for a future
theory of universal constructions for stacks.

\section{Gromov-Witten theory of varieties and stacks}
\subsection{Quantum intersection theory} In standard treatments of
Gromov-Witten theory in symplectic geometry, 
as well as \cite{Chen-Ruan}, one defines the quantum product in terms of
``$n$-point Gromov-Witten numbers'', using the inverse matrix of the
``2-point product'' on cohomology, both for defining the product and
proving associativity. In the context of algebraic cycles
on varieties, let alone stacks, this cannot be done in general, since
there is no analogue of Poincar\'e duality. A formalism which
circumvents this issue was described in the paper
\cite{Graber-Pandharipande}: 

Let $X$ be a smooth projective variety. Consider the correspondence

$$\begin{array}{ccc}
\ocM_{g,n+1}(X,\beta)&\stackrel{e_{n+1}}{\to}&X\\[2mm]
e_{1,\ldots,n}\dar             &                 &  \\[2mm]
\quad\quad\quad X^n &&\end{array}$$

\noindent where $e_{1,\ldots,n}:\ocM_{g,n+1}(X,\beta)$ is the evaluation map at
the first $n$ points, and $e_{n+1}$ is the last evaluation map.

Define a map 
\begin{eqnarray*}
A^*(X)^n & \lrar & A^*(X) \\
 \gamma_1\times\cdots\times\gamma_n & \mapsto & \left\langle
 \gamma_1,\ldots,\gamma_n,* \right\rangle_{g,\beta} , 
\end{eqnarray*}
where
$$
 \left\langle
 \gamma_1,\ldots,\gamma_n,* \right\rangle_{g,\beta} 
=({e}_{n+1})_*\ 
\left(\,e_{1,\ldots,n}^*(\gamma_1\times\cdots\times\gamma_n)
\,\cap\,[\ocM_{g,n+1}(X,{\beta})]^v\,\right) $$

\noindent Here we used the notation $[\ocM_{g,n+1}(X,{\beta})]^v$ for
the {\em virtual fundamental class} of $\ocM_{g,n+1}(X,{\beta})$,
defined using the cotangent complex of deformation theory, see, e.g.,
\cite{Behrend-Fantechi}, \cite{Li-Tian}.
 
Gromov--Witten invariants of smooth projective varieties satisfy a
number of relations which are coded by a number of well known
axioms. In particular, the genus 0 invariants can be combined to give
the multiplication law of the associative {\em quantum cohomology
ring}. For instance, the so-called ``small'' quantum product is
defined by 
$$ \gamma_1 * \gamma_2 = \sum_{\beta\in A^1(X)} \left\langle \gamma_1,\gamma_2,*\right\rangle_{0,\beta}\, q^\beta.$$

The definition using pullback and pushforward avoids entirely the use
of the inverse matrix of the multiplication table in cohomology. Of
course, in this setting we identify ``homology'' and ``cohomology''
cycles on the 
smooth variety $X$. This means that, for calculations in enumerative
geometry the inverse matrix would be used anyway.   

To prove the axioms of Gromov--Witten invariants, and, in particular,
to prove associativity of the quantum product, one needs to 
compare the fiber product
$$\ocM_{g_1,n_1+1}(X,{\beta_1})\ \mathop\times\limits_{X}\
\ocM_{g_2,n_2+1}(X,{\beta_2}), $$
along with its virtual fundamental class, with a divisor in 
$$\ocM_{g_1+g_2,n_1+n_2}(X,{\beta_1+\beta_2}),$$ as we will
illustrate later for $g=0$. 

In order to generalize this picture, we will need
\begin{enumerate} 
\item an analogue of $A^*(X)$ for a smooth Deligne--Mumford stack,
\item an analogue of the evaluation maps $e_{1,\ldots,n}$ and
$e_{n+1}$,
\item an analogue of the virtual fundamental class, and
\item an analogue of the divisor description of the fibered product.
\end{enumerate}

\subsection{Intersection theory of a stack} The first ingredient in the
``non-orbifold'' case which needs a revision in the stack situation is
the intersection theory $A^*(X)$.

Two important ingredients in intersection theory of smooth varieties are
\begin{enumerate}
\item a ring $A^*(X)$, functorial under pull-back along a morphism
$f:X\to Y$, and
\item push-forward maps $f:X\to Y$ for $f$ proper, satisfying the
projection formula.
\end{enumerate}

It has been observed long ago that such a perfect theory {\em with
integer coefficients} cannot exist for Deligne--Mumford stacks. The
problem lies with the pushing forward. Let $G$ be a finite group of
$r$ elements,
and let $p = \Spec k$ be a point. Consider the fiber diagram
$$\begin{array}{rcl}  
G    & \to & p    \\ 
\phi\dar &     & \dar \psi \\ 
p    & \to & \cB G.
\end{array} $$  
Since the degree of $G \to p$ is the order $r$ of $G$ (i.e. $\phi_*[G]
= r[p]$), and since $p
\to \cB G$ is proper and flat, it follows by the projection formula
that the degree of $\psi:p \to 
\cB G$ is $r$ as well.  But we
also have the map $\pi: \cB G \to P$, since $p$ is the coarse moduli
space of $\cB G$, and of course $\pi\circ \psi = id: p \to p$. So the
degree of $\pi$ must be $1/r$, so $r \pi_* [\cB G] = [p]$, which is
impossible since the generator $[p]$ is not divisible by $r$.

In \cite{Gillet} and \cite{Vistoli}, Gillet and  Vistoli defined an
analogous intersection  
theory {\em with rational coefficients} for smooth Deligne-Mumford
stacks. Basically
\begin{enumerate}
\item as a group, $A^*(\cX)_\QQ$ is isomorphic to
  $A^*(X)_\QQ$, where $X$ 
is the coarse moduli space, and it has the structure of a commutative
ring, but 
\item push forward is a bit more subtle, and takes into account the
degrees of stabilizers. 
\end{enumerate}

This construction can be generalized in the sense of
\cite{Fulton}, Chapter 17 (ses also \cite{Behrend-Fantechi})
to a bivariant intersection theory of Deligne--Mumford stacks - we
will use this when 
proving associativity of the quantum product.

In \cite{Kresch}, A.\ Kresch defined an intersection theory
$A^*(\cX)$ with
{\em integer coefficients} but where push-forward maps are defined
only for {\em representable morphisms} $f: \cX \to \cY$. We recall
that a morphism of algebraic stacks is representable if and only if
for any geometric object $\xi\in\cX(\Spec \CC)$, the map $f:\Aut\xi
\to \Aut\ f(\xi)$ is a monomorphism. We will use Kresch's theory when
we define an integer-coefficients version of Chen and Ruan's stringy
cohomology ring of $\cX$.

The constructions of Kresch and Vistoli are related by the fact that
$A^*(\cX) \otimes \QQ = A^*(X)_\QQ$.

\subsection{$A^*(X)$ is too small}
Both constructions discussed above are not quite sufficient for what we need. Even
for global quotients, it has long been observed that the $K$ theory of
$\cX = [V / G]$, namely the equivariant $K$-theory of $V$, is bigger
than the $K$ theory of $V/G$, even with rational coefficients. Indeed,
the decomposition theorem says that the equivariant $K$-ring $K_0(V
/\!/ G)$ is a product of 
$K_0(V^g/C(g))$, suitably twisted by cyclotomic rings. Thus any reasonable
Riemann-Roch formula must take into account contributions from
all the fixed-point sets. These contributions are lost in
$A^*(X)_\QQ$. A detailed theory of Riemann-Roch type theorems for
stacks, in which such contributions are introduced, was developed in
Toen's work \cite{Toen}. 

In enumerative geometry, T.\ Graber studied in \cite{Graber} the number of
hyperelliptic curves of genus $g$ and degree $d$ in $\PP^2$ through
$3d+1$ general points. It is explained there that the cohomology of
$Sym^2 \PP^2$ does not carry enough information to encode all these
invariants, and therefore the Hilbert scheme $Hilb_2(\PP^2)$ was used
instead. The extra information in $Hilb_2(\PP^2)$ is completely
accounted for by contributions of fixed-point sets. It is now known
(see \cite{B-K-R}) that the 
derived categories of $Hilb_2(\PP^2)$ and the {\em stack} symmetric
square $[(\PP^2)^2/ 
(\ZZ/2\ZZ)]$ are equivalent, and therefore their $K$ groups are
isomorphic.

In \cite{Chen-Ruan}, Chen and Ruan were motivated by string theory to
introduce similar contributions of fixed-point sets. We follow their
approach, with a slight modification which avoids any choice of roots of 1.

\subsection{The inertia stack and variants}
The stringy cohomology of $\cX$ according to Chen-Ruan is, as a vector
space,
the cohomology of another smooth stack, namely the inertia stack
\begin{eqnarray*}
\cX_1^{\mbox{Chen-Ruan}}=I_\cX &=& \left\{\ \ (x,g)\ \ \mid\ x\in \cO b(\cX),\  g\in \Aut\ x\right\}\\
&=&\left\{(x,H,g) \mid\ x\in \cO b(\cX),\ H\subset \Aut\ x, \ g \mbox{ a
generator of } H \right\} \end{eqnarray*} 

A slightly more canonical Gromov-Witten formalism arises if instead one takes a
Galois twist of this:
\begin{eqnarray*}
\cX_1
&=&\left\{(x,H,\chi) \mid\ x\in \cO b(\cX),\ H\subset \Aut\ x,\ \chi:H
\tilde\to\bmu_r \mbox{ for some } r\right\} \\ 
&=&\bigcup_r HomRep(\cB\bmu_r,\cX)
\end{eqnarray*} 
i.e., the stack of representable morphism from a constant cyclotomic
gerbe to 
$\cX$. 

Of course, over $\CC$ there is a natural isomorphism of $\bmu_r$ with
$\ZZ/r\ZZ$, which gives an isomorphism $\cX_1^{\mbox{Chen-Ruan}}\simeq
\cX_1$. But over other fields the group-schemes are 
different, in which case one has to use $\bmu_r$: with $\ZZ/r\ZZ$ one
cannot even define the Gromov-Witten invariants of $\cX$.

There is an interesting variant of this stack, which we denote by
$\ocX_1$. It arises as a {\em rigidification} of $\cX_1$. The point
is that the group $\bmu_r$ acts on every object of $\cX_1$ in the
component corresponding to $\bmu_r$, and the rigidification process (see
\cite{A-C-V}) allows one to remove this action from the picture (just
as the Picard functor is obtained by removing the $\CC^*$ action from
the stack of line bundles).

\begin{eqnarray*}
\ocX_1&=&\left\{(x,H,\chi) \mbox{  up to  }  \bmu_r\right\}\\
 &=& \mbox{  stack of cyclotomic gerbes in  } \cX
\end{eqnarray*} 

The geometric objects of $\cX_1$ and $\ocX_1$ are the same
$(x,H,\chi)$, but the automorphisms are different: in the first case an
automorphism is an element of the centralizer $C(H)$, in the second it
is an element of $C(H)/H$.

We note that 
\begin{itemize}
\item there are obvious inclusions $\cX \subset \cX_1$ and $\cX \subset
\ocX_1$. 
\item There is also a ``forgetful map'' $\cX_1 \to \cX$, but there
is in general no map from $\ocX_1$ to $\cX$. 
\item Finally the
rigidification map $\pi_1:\cX_1 \to \ocX_1$ exhibits $\cX_1$ as the
universal gerbe over $\ocX_1$. 

This implies that the coarse moduli
space $X_1$ of $\cX_1$ is also the moduli space of $\ocX_1$.
\end{itemize}


The total stringy Chow group of $\cX$ is defined as the total Chow group
$A^*(X_1)_\QQ$. We can also consider the integral version $A^*(\cX_1)$
defined by Kresch.

There is a useful endomorphism on 
$$r:A^*(X_1)_\QQ\to A^*(X_1)_\QQ,$$
induced by the locally constant function 
$$r:X_1 \to \ZZ$$ with value $r$
on the 
piece corresponding to $ HomRep(\cB\bmu_r,\cX)$.

\subsection{Evaluation maps}

Given a family of $n$-pointed twisted stable maps parametrized by a
scheme $S$: $$(\cC \to S, f: \cC \to \cX)$$
we have, for  each $i$ in the range $ 1\leq i
\leq n$, a map of the marking $(\Sigma_i)$ into $\cX$. This means that
$S$ parametrizes a family of cyclotomic gerbes in $\cX$. By definition
this gives a morphism $(\bar e_i)_S:S \to \ocX_1$, and since this is
true for every 
twisted stable map, we get a map of moduli stacks
$$\bar e_i:\cK_{g,n}(\cX,\beta) \to  \ocX_1.$$
Passing to moduli spaces we get 
$$e_i^{\mbox{\tiny coarse}}:\bK_{g,n}(\cX,\beta) \to  X_1.$$

Chen and Ruan \cite{Chen-Ruan} use this latter map $e_i^{\mbox{\tiny
coarse}}$  of coarse moduli
space. It so happens that the formalism of Gromov-Witten theory
basically goes
through  (with an important twist to be introduced
soon) if one works with this map {\em as if it
came from a map } $e_i^{\mbox{\tiny orb}}:\cK_{g,n}(\cX,\beta) \to
\cX_1$: it is a fact of 
life that such a ``virtual map'' does not exist in general. However, as long as
one  works with cohomology or intersection theory {\em with  
 rational coefficients,} then one can define a ``cohomological
pull-back'':
$$(e_i^{\mbox{\tiny orb}})^* = r^{-1} \cdot\bar e_i^*\  (\pi_1)_* .$$
(here again $r$ denotes the operator
which multiplies a class on a component of index  $r$ by the integer
$r$). A similar trick works also for pushforward
$$(e_{i}^{\mbox{\tiny orb}})_* = r\cdot \pi_1^*\ (\bar e_i)_*.$$

 One needs to be a little careful dealing with the  fibered products that arise
in the proof of associativity. 

In order to avoid the confusion of working with such ``weighted maps''
we can simply work with the map $\bar e_i:\cK_{g,n}(\cX,\beta) \to  \ocX_1$,
changing the products  by a suitable factor. For the sake of
exposition, it may be
a bit more transparent to use $\cX_1$ and work with   the following formalism:

Define 
$\ocM_{g,n}(\cX,\cdot)$ to be the stack of  twisted stable maps {\em with
sections of all the gerbes}. This is simply the fibered product of
the $n$ universal gerbes over $\cK_{g,n}(\cX,\cdot)$. 

Now we have evaluation maps 
$$\begin{array}{ccc}
\ocM_{g,n}(\cX,\cdot)&\stackrel{e_i}{\to}&\cX_1\\
\dar               &                 & \dar \\
\cK_{g,n}(\cX,\cdot) &\stackrel{\bar e_i}{\to}&\ocX_1\end{array}$$
 
There is an obvious involution $\iota$ on $\ocX_1$ and $\cX_1$ via $\chi
\mapsto \chi^{-1}$. We define the {\em twisted evaluation map} 
 $\check{e}_i = \iota\circ e_i$. The twisted evaluation map comes
about when one glues two twisted stable maps along two matkings to
form a {\em balanced} twisted stable map.

We can work with $\ocM_{g,n}(\cX,\cdot)$, as long as we remember to
account for the degree $1/(r_1\cdots r_n)$ of the map $\ocM_{g,n}(\cX,\cdot)
\to \cK_{g,n}(\cX,\cdot)$. Note that this degree is only locally
constant, and varies from one
connected component to another.

\subsection{Deformation and obstruction}
We want to define an appropriate virtual fundamental class for
$\ocM_{g,n}(\cX,\cdot)$, generalizing the situation of the space of
maps into a smooth projective variety.

The first obstacle one needs to overcome is mostly psychological -
does deformation theory as we know it work in the situation of a stack?

The answer is yes - and already in the literature. The point is, that
Illusie's work on deformation theory \cite{Illusie} works by deforming
a ringed topos (essentially the category of sheaves on a site), as
soon as it is reasonable enough to allow the definition of a cotangent
complex. The fact that a satck has such a topos is discussed in
\cite{L-MB}, Chapter 12, and the existence of a cotangent complex is discussed 
in \cite{L-MB}, Chapter 17 (in fact the case of Deligne--Mumford
stacks can be more easily deduced from the case of schemes).

So, as in the case if schemes, the infinitesimal deformations of a
twisted curve $\cC$ are measured by $\Ext^1(\Omega^1_\cC, \cO_\cC)$,
and if one wants to  deform some additional non-twisted marking, then
infinitesimal deformations are in $\Ext^1(\Omega^1_\cC(\log D),
\cO_\cC)$, where $D$ is the divisor given by these
markings. Obstructions would lie in the corresponding $\Ext^2$ group,
which turns out to be zero for twisted pointed curves (see
\cite{A-C-V}).

Similarly, infinitesimal deformations of a map $f: \cC \to \cX$ fixing
the structure of $\cC$ (as well as of $\cX$) are measured by
$\Hom(f^*\Omega^1_\cX, \cO_\cC) = H^0(\cC, f^* T_\cX)$. Obstructions
lie in $H^1$. 

Finally, these can be put together in the following standard manner:

Consider the complex
$$ \LL_f = [ f^* \Omega^1_{\cX} \to \Omega^1_{\cC}(\log D) ], $$
with the term on the right positioned in degree 0.
Then infinitesimal deformations of $f: \cC \to \cX$ (allowing $\cC$ and
the markings to deform) are measured by $${\operatorname{{\mathbb
E}xt}}^1(\LL_f, \cO_\cC).$$
Obstructions lie in the corresponding ${\operatorname{{\mathbb
E}xt}}^2$ group.

Now, the formalism of the virtual fundamental classes (see \cite{Li-Tian},
\cite{Behrend-Fantechi}, \cite{Kresch}) automatically allows one to
define a virtual 
fundamental class 
$[\cK_{g,n}(X,\cdot)]^v$. There is similarly a class
$[\ocM_{g,n}(\cX,\cdot)]^v$, which is just the pullback of
$[\cK_{g,n}(X,\cdot)]^v$. As we indicated before, to get the degrees
right we need to replace
the latter by a multiple: 
$[\ocM_{g,n}(X,\cdot)]^w = r_1\cdots r_n [\ocM_{g,n}(X,\cdot)]^v$.

Since deformations of twisted  marked curves are
unobstructed, it is known (see \cite{Behrend-Fantechi},
\cite{Behrend}, \cite{Bryan-Leung}) that this  
virtual fundamental class can also be constructed ``relatively to the
deformation space of the marked curve'', namely, in terms of
$\Ext^\bullet(f^*\Omega^1_\cX, \cO_\cC)$ only. This is useful when showing
associativity of the quantum product.

\subsection{The product} We can now define the Gromov-Witten operation
just as it is done for smooth projective varieties. Consider the diagram

$$\begin{array}{ccc}
\ocM_{g,n+1}(\cX,\cdot)&\stackrel{\check{e}_{n+1}}{\to}&\cX_1\\[2mm]
e_{1,\ldots,n}\dar             &                 &  \\[2mm]
\quad\quad\quad\cX_1^n &&\end{array}$$

Define 
\begin{eqnarray*}
A^*(\cX_1)^n & \lrar & A^*(\cX_1) \\
 \gamma_1\times\cdots\times\gamma_n & \mapsto & \left\langle
 \gamma_1,\ldots,\gamma_n,* \right\rangle_{g,\beta} , 
\end{eqnarray*}
where
$$
 \left\langle
 \gamma_1,\ldots,\gamma_n,* \right\rangle_{g,\beta} 
=\check{e}_{n+1\ *}\ 
\left(\,e_{1,\ldots,n}^*(\gamma_1\times\cdots\times\gamma_n)
\,\cap\,[\ocM_{g,n+1}(\cX,{\beta})]^w\,\right). $$

Again - the reason we are using the twisted evaluation map $\check{e}_{n+1}$
has to do with the fact that the twisted  curves are {\em balanced}, and will
become more explicit below.
\section{Associativity}

As in the now ``classical'' case of smooth varieties, the most subtle
property of the quantum product is its associativity. Here we go
through the main steps in the proof. We revisit these steps in the
next section when discussing the stringy cohomology ring of $\cX$.

\subsection{Associativity: the product diagram} \label{Sec:associativity-product}

For simplicity of the discussion, we restrict  from now on to the case
where $g=0$ and
$n=3$. 

 We want to prove that the ``3-point'' quantum product
is associative. Since it is commutative, it suffices to show

\begin{theorem}
$$\sum_{\beta_1+\beta_2=\beta} \left\langle\, \langle
\gamma_1,\gamma_2,*\rangle_{0,\beta_1},\, \gamma_3,\,*\,\right\rangle_{0,\beta_2}
= \sum_{\beta_1+\beta_2=\beta} \left\langle\, \langle
\gamma_1,\gamma_3,*\rangle_{0,\beta_1},\, \gamma_2,\,*\,\right\rangle_{0,\beta_2}.$$
\end{theorem}

 Fix two homology classes $\beta_i \in H_1(X, \ZZ)$,
and   for simplicity  write $\cM_i = \ocM_{0,3}(\cX,
\beta_i)$ for the corresponding stacks of twisted stable maps. The
corresponding evaluation maps are denoted $e_i^1:\cM_1 \to \cX_1$ and
$e_i^2:\cM_2 \to \cX_1$. 

The relevant diagram is the following:

$$\begin{array}{ccccc} 
                      &   &\cM_2
&\mathop\to\limits^{\check{e}_3^2}&\cX_1\\ 
                      &   &{^{e_1^2}}\swarrow \quad\searrow
{^{e_2^2}}&   &     \\ 
     \cM_1            &\mathop\to\limits^{\check{e}_3^1}&\cX_1
\quad\quad \quad\cX_1&   &     \\ 
{^{e_1^1}}\swarrow \quad\searrow {^{e_2^1}}&   &
&   &     \\ 
\cX_1 \quad\quad \cX_1&   &                      &   &       
\end{array}
$$

\ifthenelse{\equal{\version}{SHORT}}{As a first step we complete the top left corner of 
the diagram: }
{Denote  $c = (e_1^1)^* \gamma_1 \cap (e_2^1)^*
\gamma_2$ and $d=(e_2^2)^* \gamma_3$. The first step is to concentrate
on the top left corner of 
the diagram }
$$\begin{array}{ccc}
 \cM_1\,\mathop\times\limits_{\cX_1}\,\cM_2 &\mathop{\to}\limits^{p_2} & \cM_2 \\
p_1 \dar \quad & & \quad \dar e_1^2 \\
\cM_1 &\mto^{\check{e}^1_3} & \cX_1
\end{array} 
\quad \mbox{equivalently} \quad 
\begin{array}{ccc}  
\cM_1\,\mathop\times\limits_{\cX_1}\,\cM_2 & \stackrel{p_1\times p_2}\lrar              & \cM_1\times\cM_2 \\
\dar &                        &\quad \quad \quad \dar \check{e}^1_3 \times  e_1^2 \\
\cX_1&\mathop{\to}\limits^{\delta_{X}}  & \cX_1^2
\end{array}$$
\ifthenelse{\equal{\version}{SHORT}}{}{
 and understand the operation
\begin{eqnarray*}
A^*(\cM_1) &\to& A_*(\cM_2)\\
c & \mapsto &\left(\,(e^2_1)^*(\check{e}^1_3)_* (c\cap
[\cM_1]^w)\,\right)\ \cap \ 
([\cM_2]^w \cap d)
\end{eqnarray*}

We quote the following lemma from \cite{Graber-Pandharipande}:
\begin{lemma}\label{Lem:push-pull}
Consider a fiber diagram of algebraic stacks
$$\begin{array}{rlcrl}
      &            & U\times_XV &             &  \\
      &{}^{p_U}\swarrow &            &\searrow^{p_V} &  \\
      &U           &            &    V        & \\
      &{}_{q_U}\searrow &            &\swarrow_{q_V} &  \\
      &            &   X        &             &
\end{array}$$
and the associated pullback diagram
$$\begin{array}{ccc}  
U\times_XV & \to              & U \times V \\
\dar &                        &\dar \\
X&\mathop{\to}\limits^{\delta_{X}}  & X^2
\end{array}$$  
where $q_U$ and $q_V$ are proper. 
Let $A\in A_*(U), B\in A_*(V)$ and $c\in A^*(U), d\in A^*(V)$. 
Then 
$$ \left( q_V^* \ (q_U)_*\ (c\cap A)\right)\ \cap \ (B \ \cap \ d) = 
(p_V)_*\left(\ \delta_X^! (A\times B)\ \cap\ (p_U^*c\cup p_V^*d)\ \right)$$
\end{lemma}

We apply this lemma when $U=\cM_1, V = \cM_2, A = [\cM_1]^w, B =
[\cM_2]^w$, and $c,d,$ as above. We obtain
$$(e^2_1)^*(\check{e}^1_3)_* (c\cap [\cM_1]^w) \cap
([\cM_2]^w \cap d) = (p_2)_*\left(\ \delta_{\cX_1}^!([\cM_1]^w\times
[\cM_2]^w) \ \cap\ (p_1^*c\cup p_2^*d)\ \right).$$
}

Consider the following  morphisms from $\cM_1\times_{\cX_1} \cM_2$ to $\cX_1:$
\begin{eqnarray*}
e_1^\times &=& e_1^1 \circ p_1 \\ 
e_2^\times &=& e_2^1 \circ p_1 \\
e_3^\times &=& e_1^2 \circ p_2 \\
e_4^\times &=& e_1^2 \circ p_2
\end{eqnarray*}

\ifthenelse{\equal{\version}{SHORT}}{A formal push-pull argument
  detailed in \cite{Graber-Pandharipande} gives:}{  
We obtain} 
$$\left\langle \langle
\gamma_1,\gamma_2,*\rangle_{0,\beta_1},\gamma_3,*\right\rangle_{0,\beta_2}= 
(\check{e}_4^\times)_*\left(\  
\left(\ (e_1^\times)^*\gamma_1\, \cup\, (e_2^\times)^*\gamma_2\, \cup \, 
(e_3^\times)^*\gamma_3\ 
\right)\, \ \cap\, \    
\delta^!(\ [\cM_1]^w\times[\cM_2]^w\ )  \  \right)$$

As we will see,  the class $\delta^!(\ [\cM_1]^w\times[\cM_2]^w\ )$
has an interpretation as a weighted virtual fundamental class of
$\cM_1\times_{\cX_1} \cM_2$ with respect to a natural obstruction
theory. In analysing the weight, the following  
locally constant function becomes useful: 
$$ r_\times = p_1^* r_3^1.$$ 
By definition of the fibered product we also have
$$ r_\times = p_2^*  r_1^2.$$
Again, this fuction induces an endomorphism of $A_*(\cM_1\times_{\cX_1}
\cM_2)$, also denoted $r_\times$.

\subsection{Associativity: the divisor diagram}\label{Sec:associativity-divisor} The stack
$\cM_1\mathop\times\limits_{\cX_1}\cM_2$ can be viewed as the stack of
pairs of 3-pointed genus 0 twisted stable maps with {\em balanced
gluing data} along $\check{e}_3^1$ and $e_1^2$, and trivialization of
all the gerbes. Writing $\beta = \beta_1+\beta_2$, this means that
there is a gluing morphism 
$$gl:\bigcup_{\beta_1+\beta_2 = \beta}\cM_1\mathop\times\limits_{\cX_1}\cM_2 \to
\cM_{0,4}(X,\beta),$$ just as in the case of usual stable
maps.  

Consider the stabilization morphism $st:\cM_{0,4}(X,\beta) \to
\obM_{0,4}$ and the divisor $D := (12|34) \subset  \obM_{0,4}$
(corresponding to the reducible stable 4-pointed curve of genus 0
where the two first points are on the same component). We have the top
and bottom of the
following fiber diagram:

$$\begin{array}{ccccc}
& & & &\displaystyle\bigcup_{\beta_1+\beta_2 = \beta}\cM_1\mathop\times\limits_{\cX_1}\cM_2\\
& & & &\ \ \dar g                                     \\
\ocM_{0,4}(\cX,\beta)&\supset&D\mathop\times\limits_{\obM_{0,4}}\ocM_{0,4}(\cX)&
\supset& D(\cX) \\  
            \dar &   & \dar           &   & \dar \\
\ocM_{0,4}^{orb}&\supset&D\mathop\times\limits_{\obM_{0,4}}\ocM_{0,4}^{orb}&
\supset& D^{orb} \\  
            \dar &   & \dar           &   &  \\
\obM_{0,4} &\supset   & D           &   &
\end{array}$$

The notation $D(\cX)$ stands for the substack of reducible {\em twisted}
stable maps to $\cX$ with the first two marking separated from the last
two by a node. We also denote by $e_i:\ocM_{0,4}(\cX,\beta) \to \cX_1$ the
evaluation maps. 

On the second row we inserted the Artin stack $\ocM_{0,4}^{orb}$ of
{\em quasistable twisted 
curves,} i.e. connected proper twisted curves without a stability
requirement. This is a smooth stack which is neither separated nor of
finite type, yet it comes in handy just as the usual stack of
quasistable curves was handy in \cite{Behrend}. We denote by 
$$j:
D^{orb} \subset \ocM_{0,4}^{orb}$$
 the 
embedding of the substack of reducible quasistable curves with the
first two marking separated from the last two by a node.

Some relatively simple facts we  have are:
\begin{enumerate}
\item  the map $g:\bigcup \cM_1\mathop\times\limits_{\cX_1}\cM_2 \to
D(\cX)$ has locally constant degree $1/r_\times$. In fact
$\bigcup \cM_1\mathop\times\limits_{\cX_1}\cM_2$ is the universal gerbe over
the node of the universal curve of $D(\cX)$.
\item We have an equality of divisors on $\ocM_{0,4}^{orb}$:
$$[D\mathop\times\limits_{\obM_{0,4}}\ocM_{0,4}^{orb}] = r_\times [D^{orb}].$$
\item \label{It:e-factors} $e_i^\times = e_i\circ gl$.
\end{enumerate}

\subsection{Associativity: compatibility of virtual fundamental classes}\label{Sec:compatibility}
A more subtle issue is the following:
\begin{proposition} \label{Prop:compatibility}
$$\delta^!(\ [\cM_1]^w\times[\cM_2]^w\ ) = {r_\times}^2\ \cdot \ g^*
(j^! [\ocM_{0,4}(X)]^w)$$ 
\end{proposition}

{\em Proof.} By definition of $[\cM_1]^w,[\cM_2]^w, [\ocM_{0,4}(X)]^w,$ and
$ r_\times$ we need to show $$\delta^!(\ [\cM_1]^v\times[\cM_2]^v\ ) = \ g^*
(j^! [\ocM_{0,4}(X)]^v).$$ 
\ifthenelse{\equal{\version}{SHORT}}{
Denote $$\cY = \cM_1 \times _{\cX_1} \cM_2.$$ By the pullback property
of virtual fundamental classes (\cite{Behrend-Fantechi},  
Proposition 7.2) the right-hand side coincides with the obstruction
class $$[ \cY, E_{f}],$$ where $(\pi:\cC \to \cY, f: \cC \to \cX)$ is
the pullback of the universal family on $\ocM_{0,4}(X)$.

The left-hand side involves the class $[\cM_1]^v\times[\cM_2]^v$,
which, by the product property of virtual fundamental classes(\cite{Behrend-Fantechi},
Proposition 7.4), coincides with the obstruction class $$[\cM_1
\times \cM_2,E_{f_{12}}]$$ where $f_{12}$ is the disjoint union  of the pullback
of the universal twisted stable maps $(\pi_1:\cC_1 \to \cM_1,
f_1:\cC_1 \to \cX)$ on $\cM_1$ and  $(\pi_2:\cC_2 \to \cM_2,f_2:\cC_2 \to 
\cX)$ on $\cM_2$. 

Denote by $f':(p_1\times p_2)^*\cC_{12} \to \cX$ the pullback of the
disjoint union family.
}
{
The stack $\cY = \cM_1 \times _{\cX_1} \cM_2$ carries two universal families
of curves mapping to $\cX$: 
\begin{enumerate}
\item  the disjoint union $(\pi':\cC' \to \cY,f':\cC' \to \cX)$ of the pullback
of the universal twisted stable maps $(\pi_1:\cC_1 \to \cM_1,
f_1:\cC_1 \to \cX)$ on $\cM_1$ and  $(\pi_2:\cC_2 \to \cM_2,f_2:\cC_2 \to 
\cX)$ on $\cM_2$,  and 
\item the pullback $(\pi:\cC \to \cY, f: \cC \to \cX)$ by the morphism $gl$ of the
universal twisted stable map on $\ocM_{0,4}(X,\beta)$.
\end{enumerate}

By construction, there is a partial normalization  morphism $\nu: \cC'
\to \cC$ such that $f' = f\circ \nu$. The curve $\cC'$ is obtained by disconnecting
$\cC$ along a node which will be denoted $\Sigma$.  

We first interpret the right hand side of the equation. Let us apply
the pullback property of virtual fundamental classes
(\cite{Behrend-Fantechi}, 
Proposition 7.2) to the fiber diagram

$$\begin{array}{ccc}\ocM_{0,4}(X,\beta) & \stackrel{l}{\leftarrow} & D_{0,4}(X)\\
\dar && \dar \\
\ocM_{04}^{orb} & \stackrel{j}{\leftarrow} & D_{04}^{orb}
\end{array}$$

We obtain that 
$$j^! [\ocM_{0,4}(X,\beta)]^v = [D_{0,4}(X),l^*E_{f_{04}}]$$
where $E_{f_{04}}$ is the relative obstruction theory used to define
$[\ocM_{0,4}(X,\beta)]^v$, given by $\RR {\pi_{0,4}}_*f_{0,4}^*
T_\cX$.
Now, since $g:\cY \to D_{0,4}(X)$ is \'etale,
it also follows that 
$$g^*j^! [\ocM_{0,4}(X,\beta)]^v = [ \cY, E_{f}],$$
where $E_f$ is defined by $\RR \pi_* f^*T_\cX$.

We now interpret the left hand side of the equation.

We can view $\cM_1 \times \cM_2$ as a moduli stack of ``disconnected
twisted stable maps'' $(\cC_{12} \to \cM_1 \times \cM_2,
f_{12}:\cC_{12} \to \cX)$. As such, it carries a natural obstruction
theory $E_{12}$ relative to $\ocM_{03}^{orb}\times \ocM_{03}^{orb}$ defined
by 
$$  \RR{\pi_{12}}_* f_{12}^*T_\cX.$$

By the product property of virtual fundamental classes (\cite{Behrend-Fantechi},
Proposition 7.4), we have  the equality of virtual fundamental classes  $$[\cM_1
\times \cM_2,E_{12}] = [\cM_1]^v\times[\cM_2]^v.$$
}

The equality in the proposition can now be rewritten as 

$$[ \cY, E_{f}] = (\delta_{\cX_1})^! \ [\cM_1 \times \cM_2, E_{12}].$$

To prove this equality, consider the normalization exact sequence of sheaves on
$\cC$: 
$$0 \to f^* T_\cX \to \nu_*{f'}^* T_\cX \to (f^*T_\cX)|_{\Sigma} \to
0.$$

This induces a distinguished triangle on $\cY$:

$$\RR \pi_* f^* T_\cX \to \RR \pi'_* {f'}^* T_\cX \to \pi_*(\,(f^*T_\cX)|_{\Sigma}\,) \to
\RR \pi_* f^* T_\cX[1]$$

We use the terminology of \cite{Behrend-Fantechi}, Proposition
7.5. The following Lemma says that $E_f$ and $E_{12}$ are {\em
compatible 
over $\delta_{\cX_1}$,} which by \cite{Behrend-Fantechi}, Proposition
7.5 exactly implies that $[ \cY, E_{f}] = (\delta_{\cX_1})^! \ [\cM_1
\times \cM_2, E_{12}].$ \qed
\begin{lemma}$\pi_*\left((f^*T_\cX)|_{\Sigma}\right) \ \simeq\  p_1^*\,(\check e_3^1)^*\ T_{\cX_1}$
\end{lemma}

{\em Proof of the Lemma.}\marg{rewrite this?} 
 Over a geometric point $y$ of $\cY$, we can identify the
fiber of $\Sigma$ with $\cB\bmu_r$. The point $y$ maps to a geometric
point $x$ of $\cX$, with stabilizer $G$, and we can locally describe
$\cX$ around $x$ as $[U/G]$. The pulback $T$ of the tangent space of $\cX$
to $Y$ has a natural action of $\bmu_r$, and the fiber of
$\pi_*(f^*T_\cX)|_{\Sigma}$ is naturally the space of invariants
$T^{\bmu_r}$. This is naturally isomorphic to the tangent space of
$\cX_1$ at the point $\check e_3^1 (p_1(y))$, which is what we needed.
\qed

\subsection{Associativity: end of proof} Let us  apply the Proposition
to prove associativity. Denote  
\ifthenelse{\equal{\version}{SHORT}}
{ 
$$c= (e_1^\times)^*\gamma_1\ \cup\ (e_2^\times)^*\gamma_2\ \cup
\ (e_3^\times)^*\gamma_3 \ \ \in \ \ A^*( \ocM_{0,4}(X)).$$}
{ 
\begin{eqnarray*}
c&=& (e_1^\times)^*\gamma_1\ \cup\ (e_2^\times)^*\gamma_2\ \cup \
(e_3^\times)^*\gamma_3 \ \ \in \ \ A^*( \ocM_{0,4}(X))\\
A&=& [\ocM_{0,4}(X)]^w \ \ \in \ \ A_*( \ocM_{0,4}(X))\\
B=d&=& [D^{orb}]\in \ \ A^*(D^{orb}) \ \subset \ A_*(D^{orb}).
\end{eqnarray*}
We are going to apply Lemma \ref{Lem:push-pull} with $U =
\ocM_{0,4}(X), V = D^{orb}$ and $X = \ocM_{0,4}^{orb}$. We denote 
 by $ l: D(X) \hookrightarrow \cM_{04}(X,\beta)$ the embedding and by
$s:\ocM_{0,4}(X) \to \ocM_{0,4}^{orb}$ the forgetful map.}

First, it is
an easy exercise to show that $j^! [\ocM_{0,4}(X)]^w =
(\delta_{\ocM_{0,4}^{orb}})^!([\ocM_{0,4}(X)]^w \times
[D^{orb}])$.
\ifthenelse{\equal{\version}{SHORT}}
{A push-pull argument shows that }
{ Lemma \ref{Lem:push-pull} says that }
$$l_*\ (j^![\ocM_{0,4}(X)]^w \cap l^* c) = \left(s^*j_*[D^{orb}]\
\cap\ ([\ocM_{0,4}(X)]^w \cap c) \right).$$
We therefore obtain
\begin{eqnarray*}
\sum \left\langle \langle
\gamma_1,\gamma_2,*\rangle_{0,\beta_1},\gamma_3,*\right\rangle_{0,\beta_2}&=& (\check e_4)_* \ \ r_\times^2 \ \cdot \ \
\ l_* g_*
(g^*j^![\ocM_{0,4}(X)]^w \cap l^* c) \\
&=& (\check e_4)_* \ \ r_\times \ \cdot \   l_*\ (j^![\ocM_{0,4}(X)]^w \cap l^* c) \\
&=& (\check e_4)_* \ \ r_\times \ \cdot \ 
                   \left(s^*j_*[D^{orb}]\ \cap\ ([\ocM_{0,4}(X)]^w \cap c) \right) \\
&=& (\check e_4)_* \ \left(s^*(c_1(\cO_{\obM_{0,4}}(D)) \cap\ ([\ocM_{0,4}(X)]^w \cap c) \right)
\end{eqnarray*}

The latter expression is clearly independent of the way the markings
are grouped, which is what we need for the associativity theorem:
$$\sum \left\langle \langle
\gamma_1,\gamma_2,*\rangle_{0,\beta_1},\gamma_3,*\right\rangle_{0,\beta_2} = 
\sum \left\langle \langle
\gamma_1,\gamma_3,*\rangle_{0,\beta_1},\gamma_2,*\right\rangle_{0,\beta_2}.$$ \qed

\section{Stringy Chow rings with integral coefficients}

We now concentrate on the part of quantum Chow ring of $\cX$ involving only
$g=0, n=3$ and $\beta=0$. Following a suggestion of Y. Ruan during a
discussion at the conference, we 
call it the {\em stringy Chow ring} of $\cX$. 

\subsection{The main result}
A-priori, one gets a new ring structure on the group
$A^*_{st}(\cX)_\QQ := A^*(X_1)_\QQ$. It is easy to see that this is
well defined whenever $\cX$ is smooth and separated, without any
properness assumptions. Our
main claim is the following:

\begin{theorem} Let $\cX$ be a smooth and  separated Deligne--Mumford
stack, not necessarily proper. There is 
a commutative, associative ring structure on Kresch's  Chow group 
$A^*_{st}(\cX):= A^*(\cX_1)$ lifting
the   stringy chow ring structure (with rational coefficients) on
$A^*(X_1)_\QQ$.

Moreover, let $j:\cX \hookrightarrow \cX_1$ be the embedding. Then the
group homomorphism $j_*:  A^*(\cX)\to A^*(\cX_1) = A^*_{st}(\cX)$ is a
ring homomorphism, in particular $j_*1$ is the identity in   $A^*_{st}(\cX)$.
\end{theorem} 

The proof of this theorem is quite a bit more subtle than the rational case.

\subsection{Refined evaluation and the product law}
First, we need to define the product. The crucial point is that in the
case $g=0, n=3$ and $\beta=0$ we have  refined evaluation maps:

\begin{lemma}
Each evaluation maps $\bar e_i:\cK_{0,3}(\cX, 0) \to \bar\cX_1$ 
admits a lifting $\e_i:\cK_{0,3}(\cX, 0) \to \cX_1$.

The morphism $\e_i$ is representable and finite.
\end{lemma}

\proof
The existence of $\e_i$ stems from the following: given a family of
twisted stable maps $(\cC \to S, \cC \to \cX)$ in $\cK_{0,3}(\cX, 0)\
(S)$, the coarse curve is just $\PP^1\times S$, and the markings can
be identified as $\{0\}\times S, \{1\}\times S$ and $\{\infty\}\times
S$. Constructing a lifting $e_i$ is tantamount to constructing an
$r_i$-th root of the normal bundle of the $i$-th marking, functorially
in $S$, which is the same as constructing an $r_i$-th root of the
tangent space of $\PP^1$ restricted to the $i$-th marking. But a line
bundle over a point obviously has an $r_i$-th root.

The morphism $\e_i$ is easily seen to be proper and with finite
fibers. Representability is an easy calculation of stabilizers\marg{To
be written!}. 
\endproof

We note that $\cK_{0,3}(\cX, 0)$ is always smooth. The obstruction
class $O_{0,3}$ defining $[\cK_{0,3}(\cX, 0)]^v$ is given by the top Chern class of the
bundle 
$$E_{0,3}=\RR^1(\pi_{0,3})_*f_{0,3}^*T_{\cX}.$$

We can now define the product as follows:

\begin{definition}
Let $\gamma_1, \gamma_2 \in A^*(\cX_1)$. Define
$$\gamma_1\smile \gamma_2  = \check{\e}_{3*}\ 
\left(\,\e_1^*(\gamma_1)\cdot \e_2^*(\gamma_2) 
\,\cdot O_{0,3} \,\right)$$
\end{definition}

As things stand, this definition depends on the choice of the liftings
$\e_i$. However it is not difficult to show that the group
homomorphisms $\e_{i*}$ and $\e_i^*$ are independent of the choices, and
therefore the product $\gamma_1\smile \gamma_2$ is independent of the
choices.

\subsection{Associativity}
Next we need to prove associativity. The discussion in section
\ref{Sec:associativity-product} goes through word for word and one
obtains (with notation analogous to that in section \ref{Sec:associativity-product})
$$
(\gamma_1\smile \gamma_2)\smile \gamma_3 = 
(\check{\e}_4^\times)_*\left(\  
\left(\ (\e_1^\times)^*\gamma_1\, (\e_2^\times)^*\gamma_2\,  
(\e_3^\times)^*\gamma_3\ 
\right)\,  
\delta^!(\ [\cK_1]^v\times[\cK_2]^v\ )  \  \right).$$

We can now follow the arguments of section
\ref{Sec:associativity-divisor}. As before, we have a gluing morphism 
$$gl: \cK_1 \times_{\cX_1} \cK_2 \stackrel{g}{\to} D (\cX) \stackrel{l}{\hookrightarrow} \cK_{0,4}(\cX,0).$$ 
This time, the morphism $g$ is a $\bmu_{r_\times}$-bundle. 

There is a little problem  with item \ref{It:e-factors} at the end of
section \ref{Sec:associativity-divisor}: the morphism $\bar e_i:
\cK_{0,4}(\cX,0) \to \bar \cX_1$ does not lift to $\cX_1$, since the
normal bundle to the $i$-th section of the universal curve over
$\obM_{0,4}$ has degree 1, and thus has no $r_i$-th root unless
$r_i=1$. However,  we restrict
 to be the inverse image $$\tilde\cK_{0,4}(\cX,0)$$ of  the open set
$\bbA^1\, =\, \obM_{0,4}\setmin (14|23)\,  \subset\obM_{0,4}$, then a
section does exist. We denote by $gl:\cK_1 \times_{\cX_1} \cK_2 \to
\tilde\cK_{0,4}(\cX,0)$ the gluing map.
So we can 
define
$$\e_i: \tilde\cK_{0,4}(\cX,0) \to \cX_1\ , \quad i=1,\ldots,4$$
and if we do so carefully then 
$$
\e_i \circ gl = \e_i^\times.$$

Continuing with the arguments, the proof of  Proposition
\ref{Prop:compatibility} 
goes through and gives
$$\delta^!(\ [\cK_1]^v\times[\cK_2]^v\ ) = \ g^*
(j^! [\cK_{0,4}(\cX,0)]^v).$$ 
We also have that $\cK_{0,4}(\cX,0) \to \cK_{0,4}^{orb}$ is flat,
therefore
$$ g^*
(j^! [\cK_{0,4}(\cX,0)]^v) = gl^*[\cK_{0,4}(\cX,0)]^v.$$

This implies that 
$$
(\gamma_1\smile \gamma_2)\smile \gamma_3 = 
\ \check{\e}_{4*}^\times\  
 gl^*\left(\ \e_1^*\gamma_1\, \e_2^*\gamma_2\,  
\e_3^*\gamma_3\ 
 [\tilde\cK_{0,4}(\cX,0)]^v  \  \right).$$

Unfortunately, $\e_4$ is neither proper nor representable, so we
cannot replace $\check\e_{4*}^\times$ by $\check\e_{4_*} gl_*$. We
remark that Chen and Ruan utilize a certain identification of the {\em coarse}
moduli space $\bK_{0,4}(\cX,0) = \obM_{0,4} \times \bX_3$, which does
not lift to the level of stacks.

Let $\tau: \tilde\cK_{0,4}(\cX,0) \to\tilde\cK_{0,4}(\cX,0)$ be the
automorphism induced by the involution $(2,3)$ on the markings, and
denote $gl_\tau = \tau\circ gl$. One
immediately sees that associativity follows from the following:
\begin{claim} 
$$ \check{\e}_{4*}^\times\    gl^* =  \check{\e}_{4*}^\times\
(\tau\circ gl)^*:A^*(\tilde\cK_{0,4}(\cX,0)) \to A^*\cX_1. $$ 
\end{claim}

{
To show this equality it is useful to refine the divisor diagram.
 We deal with one  connected component of
$\tilde\cK_{0,4}(\cX,0)$ at a time, so let $\tilde\cK$ be such a component.
We denote $\boldsymbol{ \times} := gl^{-1} \tilde\cK$ and similarly $\boldsymbol{
 \times_\tau} := gl_\tau^{-1} \tilde\cK$. One can show that the locally
 constand function
 $r_\times$ is constant on  $\boldsymbol{ \times}$ and  $\boldsymbol{
 \times_\tau}$, and we denote these constants by $r_{\boldsymbol{\times}}$ and
 $r_{\boldsymbol{\times_\tau}}$ respectively.

We have a fiber diagram

$$\begin{array}{ccccc
}
\tilde\cK & \supset & 
        \tilde D(X) &  \leftarrow  & \boldsymbol{ \times} \\
\dar &&\dar& 
       &\dar\\
\tilde\cM_{0,4}&   
\supset &\cB \bmu_{r_{\boldsymbol{\times}}} & \leftarrow & \Spec \CC \\
\end{array}$$
and a similar diagram for the $\tau$ version.

Here $\tilde\cM_{0,4}$ denotes the open stack of twisted 4-pointed curves
over $\bbA^1$ with {\em stable} 4-pointed coarse curve, with markings
of indices
$r_1,\ldots,r_4$, and nodes of index $r_{\boldsymbol{\times}}$ over
$\{0\}\subset\bbA^1$ and index  $r_{\boldsymbol{\times_\tau}}$ over
 $\{1\}\subset\bbA^1$. It is simply obtained by endowing $\bbA^1$ with
orbifold structure of index $r_{\boldsymbol{\times}}$ at
$\{0\}\subset\bbA^1$ and index  $r_{\boldsymbol{\times_\tau}}$ at
 $\{1\}\subset\bbA^1$. It thus contains two divisors $\cG_0$ and $\cG_1$
over $\{0\}$ and $\{1\}$ in $\bbA^1$, isomorphic to
$\cB {\bmu_{r_{\boldsymbol{\times}}}}$  and $\cB
{\bmu_{r_{\boldsymbol{\times_\tau}}}}$, respectively.  

The  square on the right is cartesian since  $\Spec \CC
\to \cB \bmu_{r_{\boldsymbol{\times}}}$ is the universal principal bundle.

Ignoring the middle column, we can extend the fiber diagram as follows:

$$\begin{array}{rcl}
\tilde\cK & \stackrel{gl}{\leftarrow} & \boldsymbol{ \times} \\
s\times \check\e_4\dar &                  & \dar\check\e_4^\times \\
\tilde\cM_{0,4}\times \cX_1 & \stackrel{p\times id}{\leftarrow} &\cX_1\\
\dar &                  & \dar \\
\tilde\cM_{0,4} & \stackrel{p}{\leftarrow} &\Spec \CC
\end{array}$$
(and a similar diagram with $\tau$).

Considering the projection formula for the upper square, associativity
 follows if we prove that
$$(p\times id)^*=(p_\tau\times id)^*: A^*(\tilde\cM_{0,4}\times
\cX_1) \to A^*(\cX_1)$$
is independent

Given $\xi \in A^*(\tilde\cM_{0,4}\times
\cX_1)$, we can restrict it to the open set over $\tilde\cM_{0,4}\setmin
(\cG_0\cup \cG_1) \ \simeq \ \bbA^1 \setmin \{0,1\}$:
$$\begin{array}{ccc} & & \tilde\cM_{0,4} \\
       &\nearrow & \dar\\
\bbA^1 \setmin \{0,1\}&  \hookrightarrow & \bbA^1
\end{array}$$

 The usual exact
sequence says that this restriction extends to a class on 
$ \bbA^1 \times \cX_1$, which can be pulled back  to a class $\xi'$ on   
$\tilde\cM_{0,4}\times
\cX_1$. On the one hand,  $A^*(\cX_1
\times \bbA^1) = A^*\cX_1$, and therefore $(p\times id)^*\xi'$ is
independent of the choice of $p$.  On the other hand, $\xi$ and $\xi'$
differ by a class 
supported on $\cG_0 \cup \cG_1$. It thus suffices to show that 
$$(p\times id)^* \xi = 0$$ whenever $\xi$ is supported on $(\cG_0 \cup
\cG_1) \times \cX_1$. It suffices to treat one of the two.  
Now say $\iota: \cG_0\times
\cX_1\hookrightarrow\tilde\cM_{0,4}\times \cX_1$ is the 
embedding of the gerbe over 0, with normal sheaf $N_0$. Then
$$\iota^*\iota_*\xi =  c_1(N_0)\xi.$$
 Pulling back to $\cX_1 = \Spec \CC \times \cX_1$, we
have $$(p\times id)^*\iota_*\xi 
=(p\times id)^* c_1(N_0)\xi.$$ A simple computation shows that
$$(p\times id)^*c_1(N_0)=0,$$ 
therefore $(p\times id)^*\iota_*\xi= 0$. \qed  
}
\section{Remarks on grading, expected dimension, divisors}

\subsection{Degree shifting}
When we read about virtual fundamental classes on moduli spaces of
``usual'' stable maps, we often find it stressed that the virtual
fundamental class has the ``expected dimension''. In our discussion
here we (somewhat perversely) avoided any discussion of the expected
dimension, and evidently it is not at all necessary for proving the
associativity of the quantum product. On the other hand, the existence
of such dimension and its properties are of interest.

An interesting feature of the expected dimension  of twisted stable
maps is, that it depends on the component of the moduli space in a
subtle way, determined by the type of the evaluation maps, using the
so called ``degree shifting number'' or ``age'' of the image component
in $X_1$. 

The complex $\RR \pi_{g,n,\beta\ *} f^* T\cX$ was used
for defining the relative obstruction theory, and thus the class
$[\cK_{g,n}(\cX,\beta)]^v$. There is  a na\"\i ve approximation for the
degree of the top
Chern class,  using the usual Grothendieck--Riemann--Roch
formula. However, it is well 
known that the usual Grothendieck--Riemann--Roch formula does not hold for
Deligne--Mumford stacks; it is explained in Toen's work \cite{Toen}
that one
needs to introduce terms corresponding to pullback of the sheaf to
the twisted sectors to correct the formula. 

The issue can be understood if one considers a smooth marked twisted
curve $\cC$ with a representable map $f: \cC \to \cX$. Let $V = f^*
T_X$.
 If $\pi: \cC \to C$ is the map to the coarse curve, then $\chi (\cC,
V) = \chi( C, \pi_* V)$, since $\pi_*$ is an exact functor. Also,
$\deg_C \pi_*V = \deg_\cC \pi^*\pi_* V$. Thus the failure of
Riemann--Roch is exactly the degree of the torsion sheaf $V /
\pi^*\pi_*V$ supported at the twisted markings. 

A neighborhood of a twisted marking locally looks like $[\Spec
\CC[[t]]/\bmu_{r}]$. We describe the action as follows: we let
$\bmu_r$ act on the tangent space of $\cC$ via the fundamental
character. We can choose the parameter $t$ to be an eigenvector, and
since it generates the {\em cotangent} space at $x$, the group
$\bmu_r$ acts on $t$ via the {\em   
inverse} of the
fundamental character:$$t\ \mapsto \ \zeta_r^{-1} \ t.$$ The group
$\bmu_r$ is identified via $f$ with a subgroup 
$H$ of the stabilizer $G_x$ of a point of $\cX$. Let the tangent space
$V_x$ of $\cX$ at $x$ have basis $(v_1,\ldots,v_n)$ consisting of
$\bmu_r$-eigenvectors. The degree of $V /
\pi^*\pi_*V$ is then exactly $\frac{1}{d}\sum_i^n k_i$ where 
$$k_i = \min\{ l | v_i\cdot t^l \mbox{  is invariant  }\}.$$
In other words, $\bmu_r$ acts on $v_i$ via 
$$v_i \mapsto \zeta_r^{k_i} v_i,$$
with $0\leq k_i< r$. Thus the contribution of a marking to the degree
depends only on the component of $X_1$ where it evaluates. This
contribution $\frac{1}{d}\sum_i^n k_i$ is known as the {\em degree shifting
number} (\cite{Chen-Ruan}) or the {\em age} (Miles Reid's
terminology). If the connected components of $X_1$ are 
$X_1^i$, we denote this value by  $a(X_1^i)$.

As in \cite{Chen-Ruan:orb}, one defines a grading on the group
$A^*_{st}(\cX)$ as follows:
$$A^m_{st}(\cX) = \oplus_i A^{m-a(X_1^i)}(X_1^i).$$

An argument similar to that of Chen and Ruan shows that this makes the
stringy Chow ring with integer coefficients $A^m_{st}(\cX)$ into a
$\ZZ$-graded ring.

\subsection{Forgetting a marking}
The Gromov--Witten invariants of projective varieties satisfy the
so-called ``divisor axiom'', which comes from the fact that
$\ocM_{g,n+1}(X, \beta)$ is the universal curve over $\ocM_{g,n}(X,
\beta)$. This still works for stacks when we restrict to the case where the
$(n+1)$-st marking is {\em
untwisted}.  In our setup there is no ``forgetful map'' for a twisted
marking. There is, however, a way to change the formalism which allows
for such forgetful maps. One way to do this entails
\begin{enumerate}
\item systematically breaking up $\ocM_{g,n}(X,\beta)$ according to
the types of the markings, with components $\ocM_{g,n}(X,\beta,
\tau)$, where $\tau\in \pi_0(X_1)^n$ is an n-tuple of connected
components of $X_1$ 
indicating where each evaluation map $e_i$ lands; 
\item  for each index set $I \subset \{1,\ldots,n\}$ of {\em twisted}
markings in $\ocM_{g,n}(X,\beta, \tau)$ of size $|I| = l$ introducing
a new moduli stack 
$$\ocM_{g,n-l}(X,\beta, \tau') = [\ocM_{g,n}(X,\beta, \tau)/\cS_l].$$  
The right hand side is the quotient of $\ocM_{g,n}(X,\beta, \tau)$ by
the action of the symmetric group on $l$ letters $\cS_l$ acting on the
indices of $I$. The point is that the only thing we can ``forget''
about a collection of twisted markings is their order. The notation
$\tau'$ on the left denotes the image of $\tau$ in
$\pi_0(X_1)^n/\cS_l$ obtained by forgetting the order of the indices
in $I$.
\end{enumerate}

At this point we do not know which formalism carries more
computational or theoretical advantages.


\begin{thebibliography}{MMMMM}
\bibitem[$\aleph$-C-V]{A-C-V} D. Abramovich, A. Corti and A. Vistoli,
{\em Twisted bundles and admissible covers,} preprint 2001.

\bibitem[$\aleph$-J]{A-Jarvis} D. Abramovich and T. Jarvis, {\em
Moduli of twisted spin curves}, preprint, {\tt math.AG/0104154}



\bibitem[$\aleph$-V1]{A-V:fibered-surfaces} D. Abramovich and A. Vistoli, {\em
Complete moduli for fibered surfaces}. In {\it Recent 
Progress   in Intersection Theory}, G. Ellingsrud, 
W. Fulton, A. Vistoli (eds.), Birkh\"auser, 2000.

\bibitem[$\aleph$-V2]{A-V:families} D. Abramovich and A. Vistoli, {\em
Complete moduli for families over semistable curves}, preprint,
{\tt math.AG/9811059}.

\bibitem[$\aleph$-V3]{A-V:stable-maps} D. Abramovich and A. Vistoli, {\em
Compactifying the space of stable maps}, preprint,
{\tt math.AG/9908167}.

\bibitem[Ar]{Artin}
 M. Artin, {\em Versal deformations and algebraic stacks,} Invent. Math. {\bf
27} 
(1974), 165--189. 

\bibitem[Be]{Behrend} K. Behrend, {\em Gromov-Witten invariants in
algebraic geometry.} Invent. Math. 127   
(1997), no. 3, 601--617. 

\bibitem[B-F]{Behrend-Fantechi} Behrend, K.; Fantechi, B. The
intrinsic normal cone. Invent. Math. 128 (1997), no. 
1, 45--88. 

\bibitem[B-M]{Behrend-Manin} Kai Behrend and Yuri I. Manin, {\em  Stacks of
stable maps and  Gromov-Witten invariants,} Duke Math. J. {\bf 85} (1996), 
no.~1, 1--60.


\bibitem[B-K-R]{B-K-R} Bridgeland, Tom; King, Alastair; Reid, Miles
The McKay correspondence as an 
equivalence of derived categories. J. Amer. Math. Soc. 14 (2001), no. 3, 535--554 

\bibitem[B-L]{Bryan-Leung} Bryan, Jim; Leung, Naichung Conan The
enumerative geometry of $K3$ surfaces 
and modular forms. J. Amer. Math. Soc. 13 (2000), no. 2, 371--410

\bibitem[CR 1]{Chen-Ruan:orb} W. Chen and Y. Ruan, {\em A new cohomology
    theory of orbifolds}, preprint {\tt math.AG/0004129}

\bibitem[CR 2]{Chen-Ruan} W. Chen and Y. Ruan, {\em Orbifold
Gromov-Witten theory}, preprint {\tt math.AG/0103156}

\bibitem[D-M]{Deligne-Mumford} 
P. Deligne and  D. Mumford, {\em  The irreducibility of the space of 
curves of given genus,} Inst. Hautes \'Etudes
Sci. Publ. Math. No. {\bf 36} (1969), 75--109. 

\bibitem[F-G]{Fantechi-Gottsche} B.\ Fantechi and L.\ G\"ottsche, {\em
Orbifold cohomology for global quotients}, preprint {\tt math.AG/0104207}

\bibitem[Fu]{Fulton} W. Fulton, {\em Intersection theory.} Second
edition. Ergebnisse der Mathematik und 
ihrer Grenzgebiete,  Springer-Verlag,
Berlin, 1998.

\bibitem[F-P]{Fulton-Pandharipande} W.\ Fulton  and R.\ Pandharipande,
{\em Notes on stable maps and quantum cohomology,} in {\it Algebraic
geometry---Santa Cruz 1995}, 45--96, Proc. Sympos. Pure Math., Part 2,
Amer. Math. Soc., Providence, RI, 1997.

\bibitem[Gi]{Gillet} H.\ Gillet, {\em Intersection theory on algebraic stacks
and Q-varieties,} J. Pure Appl. Algebra 34, 193-240 (1984).

\bibitem[Gr]{Graber} T.\ Graber, {\em Enumerative geometry of hyperelliptic
 plane curves.} J. Algebraic Geom. 
10 (2001), no. 4, 725--755.

\bibitem[G-H-S]{Graber-Harris-Starr} T.\ Graber, J.\ Harris, J.\
Starr {\em Families of rationally connected 
varieties} math.AG/0109220

\bibitem[G-P]{Graber-Pandharipande} T. Graber and R. Pandharipande,
{\em The quantum product and the  intermediate Jacobian}, in preparation. 

\bibitem[Il]{Illusie}
L. Illusie, {\it Complexe cotangent et d\'eformations. I, II}, Lecture Notes in
Mathematics, Vol. 239, Springer, Berlin, 
1971 \&  Lecture Notes in
Mathematics, Vol. 283, Springer, Berlin, 1972.


\bibitem[L-MB]{L-MB} G. Laumon and L. Moret-Bailly, {\em Champs
Alg\'ebriques}.  Ergebnisse der Mathematik und ihrer Grenzgebiete 39,
Springer-Verlag, 2000. 


\bibitem[L-T]{Li-Tian} J.\ Li and G.\ Tian, {\em Virtual moduli cycles and
Gromov-Witten invariants of 
algebraic varieties.} J. Amer. Math. Soc. 11 (1998), no. 1, 119--174. 

\bibitem[K-M]{Keel-Mori} S. Keel and S. Mori, {\em Quotients by groupoids,} Ann. of
Math. (2) {\bf 145} (1997), no.~1, 193--213. 

\bibitem[Kr]{Kresch} A. Kresch, {\em Cycle groups for Artin stacks},
Inventiones Mathematicae
{\bf 138} (1999), 495--536.

\bibitem[To]{Toen} B.\ Toen, {\em Th\'eor\`emes de Riemann-Roch pour les champs
de Deligne-Mumford.} 
 $K$-Theory
18 (1999), no. 1, 
33--76. 

\bibitem[Vi]{Vistoli} A. Vistoli, {\em Intersection theory on
algebraic 
stacks and on their moduli spaces.} Invent. Math. 97 (1989), no. 3,
613-670.


\end{thebibliography}
\end{document}